\documentclass[11pt,reqno]{article}
\usepackage{euscript,amstext,amsthm,amssymb}
\usepackage{amsmath}

\newcommand{\be}{\begin{equation}}
      \newcommand{\ee}{\end{equation}}
\newcommand{\ban}{\begin{eqnarray*}}
       \newcommand{\ean}{\end{eqnarray*}}
\newcommand{\ba}{\begin{eqnarray}}
       \newcommand{\ea}{\end{eqnarray}}

\font\BBb=msbm10 at 12pt
\renewcommand{\Bbb}[1]{\mbox{\BBb #1}}
\theoremstyle{plain}
\newtheorem{Thm}{Theorem}

\newtheorem{Cor}[Thm]{Corollary}
\newtheorem{Prop}[Thm]{Proposition}
\newtheorem{Lem}[Thm]{Lemma}
\newtheorem{Def}[Thm]{Definition}
\newtheorem{Rk}[Thm]{Remark}

\newcommand{\pt}{\partial}
\newcommand{\lp}{\langle}
\newcommand{\rp}{\rangle}
\newcommand{\ra}{\rightarrow}

 \renewcommand{\o}[2]{\frac{#1}{#2}}

\begin{document}
\large
%begin Topmatter

\title{Mass under the Ricci flow}
\author{ Xianzhe Dai\thanks{Partially
supported by NSF.}\ \ and Li Ma\thanks{Partially supported by the
key project 973 of the Ministry of Sciences and Technology of
China.}} \maketitle

\begin{abstract}
In this paper, we study the change of the ADM mass of an ALE space
along the Ricci flow. Thus we first show that the ALE property is
preserved under the Ricci flow. Then, we show that the mass is
invariant under the flow in dimension three (similar results hold in
higher dimension with more assumptions). A consequence of this
result is the following. Let $(M,g)$ be an ALE manifold of dimension
$n=3$. If $m(g)\neq 0$, then the Ricci flow starting at $g$ can not
have Euclidean space as its (uniform) limit.


\end{abstract}

\maketitle
\date{}
\section{Introduction}

Ricci flow is an important geometric evolution equation in
Riemannian Geometry. It was introduced by R. Hamilton in 1982 (see
\cite{H95}) and used extensively by him to prove some outstanding
results on 3-manifolds and 4-manifolds. Recently it has been used
spectacularly by G.Perelman \cite{P} to study the geometrization
conjecture on 3-manifold. The flow has also been very useful in the
study of pinching results and metric smoothing process. As a natural
geometric tool, Ricci flow should be used to study properties of
physically meaningful objects such as mass, entropy, etc. In this
paper we would like to understand the behavior of mass under the
Ricci flow.

In general relativity, isolated gravitational systems are modeled by
spacetimes that asymptotically approach Minkowski spacetime at
infinity. The spatial slices of such spacetimes are then the so
called asymptotically flat or asymptotically Euclidean (AE in short)
manifolds. That is, Riemannian manifolds $(M^n, g)$ such that $M=M_0
\cup M_{\infty}$ (for simplicity we deal only with the case of one
end; the case of multiple ends can be dealt with similarly) with
$M_0$ compact and $M_{\infty} \simeq {\Bbb R}^n - B_R(0)$ for some
$R>0$ so that in the induced Euclidean coordinates the metric
satisfies the asymptotic conditions \be \label{afm}
g_{ij}=\delta_{ij} + O(r^{-\tau}), \ \ \  \partial_k
g_{ij}=O(r^{-\tau -1}), \ \ \
\partial_k \partial_l g_{ij}=O(r^{-\tau -2}). \ee Here $\tau >0$ is
the asymptotic order and $r$ is the Euclidean distance to a base
point. The total mass (the ADM mass) of the gravitational system can
then be defined via a flux integral \cite{ADM}, \cite{LP} \be
\label{admm} m(g)= \lim_{R \ra \infty} \frac{1}{4 \omega_n}
\int_{S_R} (\partial_i g_{ij} - \pt_j g_{ii} ) * dx_j . \ee Here
$\omega_n$ denotes the volume of the $(n-1)$-sphere and $S_R$ the
Euclidean sphere with radius $R$ centered at the base point. By
\cite{BKN89}, when the scalar curvature is integrable and $\tau >
\o{n-2}{2}$, the mass $m(g)$ is well defined and independent of the
coordinates at infinity, and therefore is a metric invariant. The
famous Positive Mass Theorem, proved by Schoen-Yau \cite{SY79}
(later Witten gave an elegant spinor proof \cite{w}), says that the
mass $m(g) \geq 0$ if the scalar curvature is nonnegative (and the
manifold is spin). Moreover, $m(g)=0$ if and only if $M$ is the
Euclidean space.

There is also the notion of asymptotically {\em locally} Euclidean,
or ALE for short, manifolds. For our purpose we will use the
following characterization of the ALE property of a complete
non-compact Riemannian manifold $(M,g)$. Namely we use the {\em
curvature decay condition}
  \be \label{cdc}
|Rm|(x)=O( d(x)^{-(2+\tau)})
  \ee
for some $\tau>0$ as $d(x) \to\infty$, and the {\em volume growth
condition}
  \be \label{vgc}
  Vol(B_r, g)\geq Vr^n
  \ee
for some constant $V>0$. Here $Rm$ is the Riemannian curvature
tensor of the metric $g$, and $d(x)$ is the distance function from
the base point.

According to \cite[Theorem (1.1)]{BKN89}, (\ref{cdc}), (\ref{vgc})
imply that $(M, g)$ is asymptotically {\em locally} Euclidean.
Namely, $M=M_0 \cup M_{\infty}$ with $M_0$ compact and $M_{\infty}
\simeq \left({\Bbb R}^n - B_R(0)\right)/\Gamma$, where $\Gamma
\subset O(n)$ is a finite group acting freely on ${\Bbb R}^n -
B_R(0)$, so that the asymptotic conditions (\ref{afm}) hold. For an
ALE manifold $(M, g)$, the mass $m(g)$ can be defined by
(\ref{admm}) again, except that $S_R$ should be taken as the
distance sphere or, equivalently, the quotient of the Euclidean
sphere by $\Gamma$.

An ALE manifold $(M, g)$ is actually AE if the asymptotic volume
ratio $\mu=1$. Here
 \be \label{avr}
\mu=\lim_{r\to\infty}V(B_r,g)/\Omega_nr^n,
 \ee
where $\Omega_n$ is the volume of the unit ball in the standard
n-dimensional Euclidean space $\mathbb{R}^n$.

As we mentioned, we would like to investigate the behavior of the
mass $m(g)$ under the Ricci flow. Recall that the Ricci flow is a
family of evolving metrics $g(t)$ such that
  \be \label{rf}
\o{\partial}{\partial t} g(t)=-2Rc(g(t)),
  \ee
on $M$ with $g(0)=g$, where $Rc(g(t))$ is the Ricci tensor of the
metric $g(t)$. To make sure that the mass is well defined under this
evolving flow $g(t)$, we first need to show that the ALE property is
preserved along the flow.

Our main result is the following
\begin{Thm} \label{mt} Let $g(t)$, $0\leq t\leq T<+\infty$,
be a Ricci flow on $M$ with $(M,g(0))$ being an ALE (AE resp.)
manifold of dimension $n$. Assume that $g(t)$ has uniformly bounded
sectional curvature. Then

 A). the ALE (AE resp.) property is preserved along the flow;

 B). the integrability condition $R\in L^1$ is also preserved along
 the flow provided (\ref{mali}) and (\ref{mali2}) hold, which is the
 case if either $R=O(r^{-q})$, $q>n$ and $\tau > \o{n-4}{2}$ or
 $\tau > \max\{n-4, \o{n-2}{2}\}$.

 C). Under the conditions above, the mass $m(t)=m(g(t))$ is well defined, and
 \be \label{crom}
m'(t)=\int_{S_{r\rightarrow\infty}}R_idS^i,
 \ee
where $dS^i=\frac{1}{4 \omega_n} * dx_i$. Furthermore, if $R=O(r^{-q})$ for $q>n$
and $\tau > \o{n-2}{2}$, or $\tau>
\max\{n-3, \o{n-2}{2}\}$, the mass is
invariant under the Ricci flow. In particular, the mass is invariant
in dimension $3$ ($\tau > \o{1}{2}$).

 D). Assume that the initial metric
$g(0)=g$ satisfies the additional decay condition
 \be \label{adc}
\partial_k \partial_l\pt_p g_{ij}=O(r^{-\tau -3}),
 \ee
then the mass is invariant under the Ricci flow if $n\leq 6$ (and
$\tau > \o{n-2}{2}$) or if $\tau > \max\{n-4, \o{n-2}{2}\}$.
\end{Thm}

That the ALE property is preserved will be obtained by using the
maximum principles of Ecker-Huisken \cite{EH91} and W.X.Shi
\cite{Shi89}, Cf. also P.Li and S.T.Yau \cite{LY}. To compute the
changing rate of the mass of the evolving metric $g(t)$, the key
part is to get a decay estimate of the space derivative of the
scalar curvature function $R(x,t)$ of the metric $g(t)$ at infinity,
which is furnished by Shi's local gradient estimate \cite{Shi89}.

%We shall assume that $(M,g)$ has finite ends throughout this paper.
% Let $(M,g)$ be
%an ALE manifold of dimension $n$ with finite ends. Recall that the
%mass of the metric $g$ is defined by
%$$
%m(g)=\int{S_{r\to{\infty}}}(g_{ij,j}-g_{jj,i})dS^i
%$$
%where $dS^i=g(\nu,\partial/\partial x^idS$ with $\nu$ being the
% outer unit normal to the sphere $S_r$ as $r\to+\infty$ ( see
%\cite{B86}).

%Our second result is

\begin{Thm} \label{st} Let $g(t)$, $0\leq t <+\infty$,
be a Ricci flow on $M$ with uniformly bounded sectional curvature.
Assume further that each $g(t)$ is ALE and $g(t)$ converges
uniformly to an ALE metric $g_{\infty}$ as $t$ goes to infinity.
Then
$$
\lim_{t\ra \infty} m(g(t))=m(g_{\infty}).
$$
\end{Thm}

The notion of uniform convergence is introduced by using the space
${\cal M}_{\tau}$ of \cite{LP}, see Definition \ref{uc}.  A direct
consequence of Theorems \ref{mt} and \ref{st} is the following

\begin{Cor} Let $(M,g)$ be
an ALE manifold of dimension $n=3$ or of asymptotic order $\tau >
\max\{ n-3, \o{n-2}{2} \}$. If $m(g)\neq 0$, then the Ricci flow
starting at $g$ can not converge uniformly to a Euclidean space.
\end{Cor}

Note that the Ricci flow preserves nonnegative scalar curvature
\cite{H95}. Thus, one can look for applications of our results by
combining with the positive mass theorem of Schoen and Yau. For
example, one sees that there are no complete non-compact Riemannian
manifolds satisfying the hypothesis of the Main Theorem in
\cite{Shi89} in dimension $3$ by using the long time convergence
result of \cite{Shi89}, Theorems \ref{mt} and \ref{st}, and the
positive mass theorem, see also \cite{CZ2}, \cite{GPZ}.

%Hence, we have the following gap result.
%\begin{Cor} \label{gt} Let $(M, g)$ be a complete noncompact
%Riemannian manifold of dimension $n \geq 3$. For any $c_1, \ c_2>0$
%and $\delta >0$, there exists a constant $\epsilon=\epsilon(n, c_1,
%c_2, \delta) >0$ such that if \\
%(A).  $$
% Vol(B(o,r)\geq c_1 r^n;
%  $$
%(B).  $$
% |\cuv m|^2 \leq \epsilon R^2, \ \ \ \ \ \ 0 \leq R \leq c_2
% (d(x))^{-(2+\delta)},
% $$
%where $\cuv m=R_{ijkl} - \o{1}{n(n-1)}R (g_{ik}g_{jl} -
%g_{il}g_{jk})$; \\
%Then $M$ must be the Euclidean space with its flat metric.
%\end{Cor}

%Thus, the complete non-compact Riemannian manifolds in the Main
%Theorem of Shi \cite{Shi89} are in fact flat. Corollary \ref{gt}
%will follow from the long time convergence result of Shi
%\cite{Shi89}, Theorems \ref{mt} and \ref{st}, and the positive mass
%theorem of Schoen and Yau \cite{SY79} (see also \cite{w}).

Let us explain why Theorem \ref{mt} comes so natural. We begin by
recalling some basic facts about Ricci flow on complete non-compact
Riemannian manifold $(M,g)$ with bounded sectional curvature $K_0$.

 Let $g(t)$ be a family of the metrics evolving under the
Ricci flow on $M$ with initial data $g$, $0\leq t\leq T<+\infty$. We
shall write by $\nabla_{g(t)}$ and $R_{ijkl}(t)$ the Riemannian
connection and Riemannian curvature tensor of $g(t)$ respectively.
R.Hamilton proved in \cite{H95} that the asymptotic volume ratio
$\mu(t)=\mu(g(t))$ of (\ref{avr}) is a constant under the Ricci flow
with bounded curvature and nonnegative Ricci curvature, where
$|Rm|\to 0$ at infinity on complete non-compact Riemannian manifold.
This result tells us that if $\mu\neq 1$ at $t=0$, then the Ricci
flow can not have Euclidean space as its limit.

It is well known that Ricci flow smoothes out the metric. W.X.Shi
\cite{Shi89-2} showed that there exists a positive constant $T>0$
such that for any integer $\alpha\geq 0$ and any $0<t\leq T$, there
exist constants $c(n,K_0)$, $c(n,K_0,T)$ and $c(n,,K_0,\alpha,t)$
such that
$$
e^{-c(n,K_0)t}g\leq g(t)\leq
e^{c(n,K_0)t}g,\;\;|\nabla_g-\nabla_{g(t)}|\leq c(n,K_0)t,
$$
$$
|R_{ijkl}|\leq c(n,K_0,T),
$$
and
$$
|\nabla^{\alpha}_{g(t)}R_{ijkl}(t)|\leq c(n,K_0,\alpha,t).
$$
All these facts will be implicitly used in this paper. We note that
there is also uniqueness in this setting by the recent result of
\cite{CZ} and \cite{DM2}. It is clear that the volume growth
condition (\ref{vgc}) is preserved along the Ricci flow. In section
18 in \cite{H95}, R.Hamilton further proved that if the curvature
$Rm\to 0$ as $s\to +\infty$ for the initial metric, where $s$ the
distance function to a fixed point of the metric $g$, the same is
true for each $g(t)$. So it is very natural for one to expect that
if the curvature of the initial metric has decay at infinity, then
the same is true for the evolving metric $g(t)$. With this
understanding, we want to know the change of the mass under the
flow.

Throughout this paper we will denote by $C, c$ various constants
depending only on dimension.

\section{Preliminaries}

In this section we briefly introduce some facts on Ricci flow. we
shall use notations from \cite{H93}. Let $M$ be a manifold of
dimension $n$, $g(t)$ a family of metrics evolving by Ricci flow
(\ref{rf}). The curvature tensor evolves by the equation
 \be \label{ctee}
 \o{\pt}{\pt t} Rm = \Delta Rm + Rm * Rm,
 \ee
where $Rm * Rm$ denotes a quadratic expression of the curvature
tensor. It follows then
 $$
\o{\pt}{\pt t} |Rm|^2 = \Delta |Rm|^2 - 2 |\nabla Rm|^2 + Rm * Rm *
Rm,
 $$
which yields
 \be \label{noctee}
 \o{\pt}{\pt t} |Rm|^2 \leq \Delta |Rm|^2 + C |Rm|^3.
 \ee

The evolution equation for the scalar curvature is much simpler, and
one has
 \be \label{scee}
 \o{\pt}{\pt t} R = \Delta R + 2|Rc|^2.
 \ee

Now let $X$ be a point in $M$. Let $Y=\{Y_a\}, 1\leq a\leq n$ be a
frame at $X$. In local coordinates $X=\{x^i\}$, we have
$$
Y_a=y_a^i\partial/\partial x^i.
$$

Let $$g_{ab}=g(Y_a,Y_b),$$ and let
$$
\nabla_b^a=y_b^i\partial/\partial y_a^i
$$
be the vector fields tangent to the fibers of the frame bundle.
Write by $D_a$ the vector field on the frame bundle $F(M)$ which is
the lift of the vector $Y_a$ at $Y\in F(M)$. Then we have
$$
D_a=y^i_a[\partial x^i-\Gamma^k_{ij}y^j_b\partial/\partial y_b^k]
$$
where $\Gamma^k_{ij}$'s are the Christoffel symbols of the
connection. Under the Ricci flow, we can define the evolving
orthonormal frame on $M$ such that
$$
\partial_tF^i_a=g^{ij}R_{jk}F^k_a,
$$
where $(g^{ij})$ is the inverse matrix of $(g_{ij})$. Then we set
$$
D_t=\partial_t+R_{ab}g^{bc}\nabla_c^b.
$$
Note that
$$
D_tg_{ab}=0.
$$
This says that $D_t$ is the unique tangent vector field to the
orthonormal bundle. Choose a metric on $F(M)$ such that $D_a$,
$\nabla_c^b$ are an orthonormal basis. Then we can see that
$D_t-\partial_t$ is a space-like vector orthonormal to the
orthonormal frame bundle. A useful fact for us is that for a smooth
function $u$ on $M\times (0,T)$, we have
 \be \label{fofb}
(D_t-\Delta)D_au=D_a(\partial_t-\Delta)u.
 \ee

  We now recall R.Hamilton's argument \cite{H95}.
  Assume that a $K$-bounded smooth function $u$ satisfies the heat equation
  \be \label{he}
u_t=\Delta u,\;\;{in}\;\; M
  \ee
  with $|Du|^2\leq \delta$ at $t=0$.
Then we have by (\ref{fofb})
$$
D_tD_au=\Delta D_au
$$
and thus
$$
\partial_t|Du|^2=\Delta |Du|^2-2|D^2u|^2.
$$
By the maximum principle of Shi \cite{Shi89} we have
\begin{equation}
|Du|^2(x,t)\leq \delta. \label{ccc}
\end{equation}
Let $F=t|D^2u|^2+|Du|^2$. Then by a direct computation, we have
$$
\partial_tF\leq \Delta F-(1-cKt)|D^2u|^2.
$$
Hence, by the maximum principle of Shi \cite{Shi89} again we get for
$t\leq 1/cK$,
$$
F(x,t)\leq \delta^2,
$$
which implies that
$$
|D^2u|\leq \delta/\sqrt{t}.
$$
Note that $|\Delta u|^2\leq n|D^2u|^2$. Using the heat equation we
obtain that for $t\leq 1/cK$,
$$
|u_t|\leq \sqrt{n}\delta/\sqrt{t}.
$$
Therefore, we have
$$
|u(x,t)-u(x,0)|\leq 2\sqrt{n}\delta\sqrt{t}.
$$
Take $\delta\leq \sqrt{K}\epsilon^2$ and assume that
$$
\lim_{s\to\infty}u(x,0)=0,\;\;\;in\;\;M
$$
uniformly. Then we can conclude using an iteration argument that for
any $t\in[0,T]$,
$$
\lim_{s\to\infty}u(x,t)=0,\;\;\;in\;\;M
$$
uniformly. In fact, R.Hamilton \cite{H95} has showed that for any
$\delta>0$ and for any bounded smooth function $u_0\in C^1(M)$ with
$\lim_{s\to\infty}u_0(x)=0$,
 one can find a bounded smooth solution $u(x,t)$ to the heat equation
 such that $u_0(x)\leq u(x,0)$ and $|Du|^2(x,t)\leq \delta$ on $M\times[0,T]$.

\section{ALE is preserved}

In this section we study the ALE property under the Ricci flow and
show that it is preserved. It can be reduced to studying the
non-negative solutions to the heat equation
\begin{equation}
u_t=\Delta u,\;\;{in}\;\; M\label{mm}
\end{equation}
with initial data $u(0)=u_0$, where $\Delta=\Delta_{g(t)}$ is the
Laplacian operator of the family of metrics $g(t)$. We assume that
$u_0$ has a decay $O(d(x)^{-\sigma})$ for some $\sigma>0$.

\begin{Thm} \label{drip} Let $g(t)$ be the solution of the Ricci flow (\ref{rf}) over
$[0, T]$. Assume that $g(t)$ has uniform curvature bound
$|Rm(g(t))|\leq K$. Then non-negative solutions to (\ref{mm}) have
the same decay rate as the initial data $u_0$.
\end{Thm}

The main tools here are the maximum principles, especially the
maximum principle of \cite[Theorem 4.3]{EH91}. For reader's
convenience, we quote the result here (the superscript `$t$' is put
in here to emphasize the $t$-dependence from the metric $g(t)$).

\begin{Thm}[Ecker-Huisken]\label{xianzhe} Suppose that the complete non-compact
manifold $M^n$ with Riemannian metric $g(t)$ satisfies the uniform
volume growth condition \be \label{uvgc} {\rm vol}^t(B^t_r(p)) \leq
\exp \left(k(1+r^2)\right) \ee for some point $p\in M$ and a uniform
constant $k>0$ for all $t\in [0, T]$. Let $w$ be a function on $M
\times [0, T]$ which is smooth on $M \times (0, T]$ and continuous
on $M \times [0, T]$. Assume that $w$ and $g(t)$ satisfy

i). the differential inequality \be \label{dihe}
\o{\partial}{\partial t}w - \Delta^t w \leq {\bf a}\cdot \nabla w +
b w, \ee where the vector field ${\bf a}$ and the function $b$ are
uniformly bounded \be \label{uboc} {\rm sup}_{M\times [0, T]}
|\,{\bf a}| \leq \alpha_1, \ \ \  {\rm sup}_{M\times [0, T]} |\,b|
\leq \alpha_2 \ee for some constants $\alpha_1, \alpha_2 < \infty$;

ii). the initial data \be \label{inid} w(p, 0) \leq 0 \ee for all
$p\in M$;

iii). the growth condition \be \label{gcos} \int_0^T \left( \int_M
\exp \left[ - \alpha_3 d^t(p, y)^2\right] |\nabla w |^2 (y) d\mu_t
\right) dt < \infty \ee for some constant $\alpha_3 > 0$;

iv). bounded variation condition in metrics \be \label{bvim} {\rm
sup}_{M\times [0, T]} |\o{\partial}{\partial t} g(t)| \leq \alpha_4
\ee for some constant $\alpha_4 < \infty$. \\
 Then, we have
 \be w \leq 0 \ee
 on $M\times [0, T]$.
\end{Thm}

In our situation, with the metric $g(t)$ coming from the Ricci flow
(\ref{rf}), the condition (\ref{bvim}) is clearly satisfied by the
uniform curvature bound. The uniform volume growth condition
(\ref{uvgc}) also follows immediately from the volume comparison
theorem via the curvature bound. The differential inequality will be
coming from a modification of the solution of the heat equation
(\ref{mm}). To see that the coefficients are uniformly bounded per
(\ref{uboc}) requires the following lemma.

\begin{Lem} \label{ml} Let $g(t)$ be the solution of the Ricci flow (\ref{rf}) over
$[0, T]$ with $g(0)=g$ being an ALE. Assume that $g(t)$ has uniform
curvature bound $|Rm(g(t))|\leq K$. Then, for sufficiently large
$R$, there is a smooth positive function $f$ on $M$ such that
$$
f(x)=C_0 > > 1, \ \ \ \mbox{for} \ x \in B_R;
$$
$$ c\,d^t(x) \leq f(x) \leq Cd^t(x) \ \ \ \mbox{for} \ x \in M- B_R.
$$ Moreover
$$ f\geq C_0, \ \ \ | \nabla^t f | \leq C_1, \ \ \ | \Delta^t f | \leq C_2
. $$
\end{Lem}

\begin{proof} Since $(M, g)$ is ALE, we have coordinates at infinity,
% \cite{BKN89}, therefore, according to \cite{B86}, we can choose it to harmonic
%coordinates and
which we denote by $x$. Let $|x|$ be the Euclidean distance
function.
%Then, for $|x| \geq 1$ (say),
% \be \label{defd}
% |\nabla_g |x|\,|\leq c_1, \ \ \ \ \ \  | \Delta_g |x|\,| \leq c_2.
% \ee

Choose a smooth increasing function $\phi(s)$ on $\mathbb R$ such
that
 \begin{align*}
\phi(s)= & C_0 =R- 1, \ \ \ \mbox{if} \ s\leq R-1; \\
\phi(s)= & s,\ \ \ \mbox{if} \ s\geq R,
 \end{align*}
and
 $$
 |\phi'| \leq 1, \ \ \ \ \ \ |\phi''| \leq 2.
 $$

We define our function $f$ to be $f(x)=\phi(|x|\,)$. Then clearly
$f\geq C_0$. Since the metrics $g(t)$ are all equivalent and
$g(0)=g$ is ALE, we can use the Euclidean norm in estimating
$|\nabla^t f|$ and $|\Delta^t f |$. Then
 $$
 |\nabla^t f | = |\phi'| |\nabla^t |x|\, | \leq C_2.
 $$
Similarly the estimate
 $$
  | \Delta^t f | \leq C_2
 $$
follows from the coordinate expression of the Laplacian
$$
\Delta^t =\o{1}{\sqrt{\det g(t)} } \o{\partial}{\partial x_i}
\left(\sqrt{\det g(t)} g^{ij}(t) \o{\partial}{\partial x_j} \right)
$$
and the known estimate for $g(t)$.
\end{proof}

We now suppress the superscript `$t$' with the understanding that
all covariant derivatives and Laplacian are taken with respect to
$g(t)$.

In our application of Theorem \ref{xianzhe} to the proof of Theorem
\ref{drip}, we will let $w=f^{\sigma}u$, where $\sigma>0$. Then the
growth condition (\ref{gcos}) follows from the gradient bound
(\ref{ccc}), which implies that
$$
|\nabla w|\leq C(T)f^{\sigma+1}.
$$

We now turn to the proof of Theorem \ref{drip}.

\begin{proof}

  For simplicity, we assume that $(M,g_0)$ is an ALE with one end.
%  Then we can find a compact
%  subset $K\subset M$ and global coordinates $(x^i)$ on $M-K$ such that the
%  decay for the initial data $u_0$ is that
Let $u_0(x)=O(d(x)^{-\sigma})$ as $d(x)\to\infty$,
% with $2<\sigma<3$.
where the distance function is with respect to a fixed point $o$ in
$M$. Choose a global smooth positive function $f(x)$ on $M$
%such that
%$$
%f(x)=C_0>>1,\;\;\;in\;\;K\subset\{d_{g_0}(x,o)<d_0<+\infty\}
%$$
as in Lemma \ref{ml} and let
$$
h(x)=f(x)^{\sigma}.
$$
%in $\{d_{g_0}(x,o)>d_0+d\}$, where $C_0,d_0, d>0$ are positive
%constants.

Set $w(x,t)=h(x)u(x,t)$. Then, by a direct computation we have that
$$
w_t=hu_t,
$$
$$ w_i=h_iu+hu_i
$$
and
$$
\Delta w=\Delta h u+2\nabla h\nabla u+h\Delta u.
$$
Hence,
$$
(\partial_t-\Delta)w= Bw-2\nabla \log h \nabla w.
$$
where $B(x,t)=\frac{2|\nabla h|^2-h \Delta h}{h^2}$. Note that the
coefficients $B$ and $\nabla \log h$ are uniformly bounded by Lemma
\ref{ml}. In particular, $|B| \leq b$.
%we have for $|x|>>1$,
%$$
%w_t-\Delta w\leq -\frac{2\sigma}{|x|}x\cdot\nabla w
%$$
%on $M_T:=M\times [0,T]$.
Since
$$w(x,0)=d(x)^{\sigma} u_0(x)\leq D<+\infty,$$
and
$$
(\partial_t-\Delta)(w-De^{tb})\leq B(w-De^{tb})-2\nabla \log h\nabla
(w-De^{tb}),
$$
%on the set $\{w(x,t)>D+\epsilon\}$ for every small $\epsilon>0$,
we have by the maximum principle of Ecker-Huisken, Theorem
\ref{xianzhe}, (see also the proofs of Theorem 18.2 \cite{H95} and
Theorem 4.3 in \cite{EH91}, see also \cite{Shi89}) that there exists
a uniform constant $C_1>0$ such that
$$
max_{M_T}w\leq C_1,
$$
where $M_T=M\times [0, T]$. This implies the desired decay for
$u(x,t)$.

\end{proof}

We are now in a position to prove the first part of Theorem
\ref{mt}. In fact, since $|Rm(g_0)|\leq C_0d(x)^{-\sigma}$ for
$d(x)>>1$, where $\sigma=2+\tau$, we can choose a bounded smooth
function $u_0$, which dominates the function $|Rm(g_0)|^2$ such that
it is $C_0d(x)^{-2\sigma}$ for $d(x)>>1$ and has bounded gradient.
Let $u$ be the solution of heat equation as above. Then
 under Ricci flow, we have from (\ref{noctee}) and the uniform
 curvature bound
$$
\partial_t |Rm|^2\leq \Delta |Rm|^2+CK|Rm|^2
$$
while
$$
\partial_t (e^{CKt}u)=\Delta (e^{CMt}u)+CKe^{CKt}u.
$$
Therefore, we have by the maximum principle of Shi \cite{Shi89} and
Theorem \ref{drip} that
$$
|Rm|^2\leq e^{CKt}u\leq e^{CKt}d(x)^{-2\sigma}, \;\;{on}\;\; M.
$$
 Thus, under the Ricci flow, the ALE property is preserved.

\begin{Rk} That the ALE property is preserved does not follow from
\cite[Remark 0.9]{K}, as was claimed there.
\end{Rk}

The same is true for AE, as we have the following analog of a
theorem of Hamilton \cite{H95}.

\begin{Cor} \label{avrt} Let $g(t)$, $0\leq t\leq T$,
be a Ricci flow on $M$ with uniformly bounded sectional curvature.
Assume further that $g(0)$ is ALE. Then the asymptotic volume ratio
$\mu(t)=\mu(g(t))$ is constant along the Ricci flow.
\end{Cor}

\begin{proof} If $(M, g)$ is ALE, it follows from the characterization
in \cite{BKN89} that $M=M_0 \cup M_{\infty}$ with $M_0$ compact and
$M_{\infty} \simeq \left({\Bbb R}^n - B_R(0)\right)/\Gamma$, where
$\Gamma \subset O(n)$ is a finite group acting freely on ${\Bbb R}^n
- B_R(0)$, so that the asymptotic conditions (\ref{afm}) hold, where
the asymptotic coordinates comes from the projection of the
Euclidean coordinates under ${\Bbb R}^n - B_R(0) \ra \left({\Bbb
R}^n - B_R(0)\right)/\Gamma$. Therefore, $(M, g)$ has the asymptotic
volume ratio
 $$
 \mu=\o{1}{|\Gamma|}.
 $$
 Since ALE is preserved along Ricci flow by Theorem \ref{mt}, we
deduce that $\mu(t)$ is a constant by the continuity of $\mu(t)$, as
shown in the following lemma.
\end{proof}

\begin{Lem} \label{cavr} Let $g(t)$, $0\leq t\leq T$,
be a Ricci flow on $M$ with uniformly bounded scalar curvature.
Assume that the asymptotic volume ratio $\mu(t)=\mu(g(t))$ are well
defined. Then $\mu(t)$ is a continuous function of $t$.
\end{Lem}

\begin{proof} The volume of a ball changes according to the formula
 \be \label{rcfvob}
 \o{d}{dt} V(B_r, g(t)) = - \int_{B_r} R(g(t)) dv_{g(t)}.
 \ee
Hence
 $$
 |\o{d}{dt} V(B_r, g(t))| \leq K V(B_r, g(t)),
 $$
where $K$ denotes the uniform bound on the scalar curvature.
Therefore,
 $$
 e^{-K(t-t_0)} V(B_r, g(t_0)) \leq V(B_r, g(t)) \leq e^{K(t-t_0)} V(B_r,
 g(t_0)),
 $$
from which the continuity follows.

\end{proof}

\section{The changing rate of mass}

For the mass to be well defined, one needs the integrability
condition $R\in L^1$ in addition to the requirement that the
asymptotic order $\tau > \o{n-2}{2}$ \cite{B86}, \cite{S}. We have
seen that the ALE property is preserved along the Ricci flow. We now
examine the integrability condition.

One thing that in particular guarantees the integrability is the
decay condition
 \be \label{scdc}
 R=O(r^{-q}), \ \ \ \ \ \ q> n.
 \ee
We have

\begin{Thm} \label{scdt}
Let $g(t)$ be the solution of the Ricci flow (\ref{rf}) over $[0,
T]$ with uniformly bounded curvature. Assume that $g(0)$ is ALE with
asymptotic order $\tau>0$ and its scalar curvature satisfy the decay
condition $R(0)=O(r^{-q}), \ \ q>0$. Then the scalar curvature of
$g(t)$ satisfies the decay condition
$$
R(t)=O(r^{-q'}), \ \ \ \ \ \ q'=\min\{ q, 2(\tau +2) \}.
$$
\end{Thm}

\begin{proof} This is similar to our proof of the ALE property. The
scalar curvature satisfies the evolution equation
 $$
 \o{\pt}{\pt t} R = \Delta R + 2 |Rc|^2.
 $$
By the assumption and our result on the ALE property, we have
$|Rc|^2=O(r^{-2(\tau +2)})$. Let $f$ be the function in Lemma
\ref{ml} and $w=f^{q'} R$. Then
 $$
 (\pt_t - \Delta)w \leq Bw + 2q' \nabla \log f \nabla w + C,
 $$
where $f^{q'} |Rc|^2 \leq C$. Hence, by the argument in the proof of
ALE property, we have
 $$
 \max_{M_T} w \leq C_1.
 $$
\end{proof}

In particular, when the order of decay of the initial scalar
curvature $q>n$ and the asymptotic order $\tau> \o{n-4}{2}$, then
the order of decay of the evolving scalar curvature also satisfies
$q'>n$.

In general, without assuming the ALE conditions, we show that under the natural condition
\begin{equation}
\int_0^T\int_M|Rc|^2<+\infty, \label{mali}
\end{equation}
the property $R\in L^1$ is preserved under the Ricci flow if
%uniformly bounded curvature.
% $R\geq 0$ and
the decay condition
 \be \label{mali2}
 R(t)=O(r^{-\sigma})
 \ee
holds uniformly for some $\sigma\geq n-2$ and all $t\in [0, T]$. We remark that both
conditions (\ref{mali}) and (\ref{mali2}) are always true for
$0<T<+\infty$ if the initial metric is ALE, provided the order $\sigma=2+\tau\geq n-2$  (i.e. $\tau \geq n-4$)
and $2\sigma>n$ in the curvature decay condition
(\ref{cdc}), as it follows from our result on ALE property and
Theorem \ref{scdt}.
%We may assume that $n\geq 4$. Then $\sigma=2+\tau>n-2\geq n/2$.

\begin{Thm} \label{isct}
Let $g(t)$ be the solution of the Ricci flow (\ref{rf}) over $[0,
T]$ with uniformly bounded curvature. Assume that the conditions
(\ref{mali}) and (\ref{mali2}) hold.
Then the property $R\in L^1$ is preserved under the Ricci flow.
\end{Thm}

\begin{proof}

Recall that on $M$,
$$
\Delta R=R_t-2|Rc|^2
$$
Let $p=1+\epsilon$ with small $\epsilon>0$. Let $\phi$ be a
non-negative cut-off function such that $0\leq \phi\leq 1$ on $M$,
$\phi=1$ on $B_r(o)$, $\phi=0$ outside $B_{2r}(o)$, and
$$
|\nabla \phi|^2\leq 4\phi/r^2.
$$
Then
\begin{align*}
&\int_0^tdt\int_M\phi^2|R|^{p-1}(R_t-2|Rc|^2)\,{\rm sgn}R\\
& = \int_0^tdt\int_M\phi^2|R|^{p-1}(\Delta R)\, {\rm sgn}R\\
&=-2\int_0^tdt\int_M\phi |R|^{p-1}<\nabla\phi,\nabla R> {\rm sgn}R
\\
& -(p-1)\int_0^tdt\int_M\phi^2|R|^{p-2}|\nabla R|^2\\
&\leq\frac{2}{p-1}\int_0^tdt\int_M|\nabla
\phi|^2|R|^p-\frac{2(p-1)}{p^2}\int_0^tdt\int_M\phi^2|\nabla(|R|^{p/2})|^2\\
&\leq \frac{2}{(p-1)r^2}\int_0^tdt\int_M\phi |R|^p
-\frac{2(p-1)}{p^2}\int_0^tdt\int_M\phi^2|\nabla(|R|^{p/2})|^2\\
&\leq \frac{2}{(p-1)r^2}\int_0^tdt\int_M\phi |R|^p
\leq\frac{2CTr^{n-2-p\sigma}}{p-1}\to 0
\end{align*}
as $r\to \infty$. Here we have used the decay condition
$R=O(r^{-\sigma})$ for some $\sigma>n-2$.

 By direct computation, we have
\begin{align*}
&\int_0^tdt\int_M\phi^2|R|^{p-1}(R_t-2|Rc|^2)\,{\rm sgn}R\\
&=-2\int_0^tdt\int_M\phi^2|R|^{p-1}|Rc|^2 \,{\rm sgn}R +
\int_0^tdt\int_M\phi^2|R|^{p-1}R_t \,{\rm sgn}R \\
&=-2\int_0^tdt\int_M\phi^2|R|^{p-1}|Rc|^2 \,{\rm sgn}R
+\frac{1}{p}\int_0^tdt\frac{d}{dt}\int_M
\phi^2|R|^p \\
& +\frac{1}{p}\int_0^tdt\int_M \phi^2|R|^{p}R\\
&=-2\int_0^tdt\int_M\phi^2|R|^{p-1}|Rc|^2\,{\rm sgn}R
+\frac{1}{p}\int_M \phi^2|R|^p(t)
-\frac{1}{p}\int_M \phi^2|R|^p(0)\\
&+\frac{1}{p}\int_0^tdt\int_M \phi^2|R|^{p}R.
\end{align*}
Hence, we have
\begin{align*}
&-2\int_0^tdt\int_M\phi^2|R|^{p-1}|Rc|^2 \,{\rm sgn}R
+\frac{1}{p}\int_M
\phi^2|R|^p(t)\\
&-\frac{1}{p}\int_M \phi^2|R|^p(0) +\frac{1}{p}\int_0^tdt\int_M
\phi^2|R|^{p}R\\
&\leq\frac{2CTr^{n-2-p\sigma}}{p-1}.
\end{align*}
 Sending $r\to +\infty$, we get
\begin{align*}
& -2\int_0^tdt\int_M|R|^{p-1}|Rc|^2 \,{\rm sgn}R +\frac{1}{p}\int_M
|R|^p(t) -\frac{1}{p}\int_M |R|^p(0) \\
& +\frac{1}{p}\int_0^tdt\int_M |R|^{p}R \leq 0.
\end{align*}
Sending $p\to 1$, we have that
\begin{align*}
-2\int_0^tdt\int_M|Rc|^2\,{\rm sgn}R +\int_M |R|(t) -\int_M
|R|(0)+\int_0^tdt\int_M |R|R \leq 0.
\end{align*}
That is,
$$
\int_M |R(t)| -\int_M |R(0)|\leq \int_0^Tdt\int_M R^{2} +
2\int_0^Tdt\int_M|Rc|^2,
$$
which implies that $R\in L^1$ for each $t>0$.
\end{proof}

Using similar argument, we can show that the property $|Rm|\in L^p$
($p\geq 1$) is preserved under the Ricci flow with bounded
curvature.

%Therefore, the ADM mass is well-defined under Ricci flow with
%initial ALE decay and $|Rm|\in L^2$, and with non-negative scalar
%curvature.

We now look at the change of mass under the Ricci flow. Let $S$ be a
hypersurface in $M$. Without loss of generality, we can assume that
$M$ is oriented.

We can take the local frame ${F_{a}}$ such that $F_{1},...,F_{n-1}$
are tangent to $S$ and $F_{n}=\nu$ is orthogonal to $S$ at $X$.

Let ${\omega^{a}}$ be a local frame dual to ${F_{a}}$. Then the area
form on $S$ is
$$
dS=\omega^{1}\wedge...\wedge\omega^{n-1}.
$$
Let
$$
f^{a}_{b}=\partial/\partial_{t}\omega^{a}(F_{b}).
$$
Then we have
$$
f^{a}_{b}=-\omega^{a}(\partial/\partial_{t}F_{b}),
$$
which is a decay term of the same order as $R_{jk}=O(r^{-\sigma})$,
and
$$
\partial/\partial_{t}\omega^{a}=f^{a}_{b}\omega^{b}.
$$
Hence,
$$
\partial/\partial_{t}dS=O(r^{-\sigma}).
$$
It is also clear that
$$
\partial/\partial_{t}g(F_{n},F_{a})=O(r^{-\sigma}).
$$
Hence
$$
\partial/\partial_{t}dS^{i}=O(r^{-\sigma}).
$$
Therefore,
\begin{align*}
m'(t)
&=\int_{S_{r\rightarrow\infty}}(\frac{\partial}{\partial
t}g_{ij,j}-\frac{\partial}{\partial t}g_{jj,i})dS^{i}\\
&+\int_{S_{r\rightarrow\infty}}(g_{ij,j}-g_{jj,i})\frac{\partial}{\partial
t}dS^{i}
\end{align*}
and the second term in the above equation is zero. So, under the
flow, we have
$$
m'(t)=-2\int_{S_{r\rightarrow\infty}}(R_{ij,j}-R_{jj,i})dS^i.
$$
By using the contracted first Bianchi identity
$$
2R_{ij,j}=R_i,
$$
we have that
\begin{align*}
m'(t)&=-2\int_{S_{r\rightarrow\infty}}(\frac{1}{2}R_{i}-R_{i})dS^i\\
&=\int_{S_{r\rightarrow\infty}}R_idS^i.
\end{align*}

Note that, by using the local gradient estimate of Shi (see Theorem
13.1 in \cite{H95}), we have that $R_i=O(|x|^{-\sigma})$ for
$\sigma=2+\tau$. Hence, we have
$$
m'(t)=0
$$
when $n=3$. The same is true if $\tau > n-3$ for any dimension
$n\geq 3$.

Similarly, if the initial metric satisfies the additional decay
condition \ref{adc}, then one has a better estimate
$R_i=O(|x|^{-\sigma})$ for $\sigma=3+\tau$. Therefore, $ m'(t)=0 $
provided $\tau+3>n-1$. This will be the case if $n \leq 6$ and $\tau
> \o{n-2}{2}$, or $\tau > \min \{ n-4, \o{n-2}{2} \}$.

On the other hand, if $R=O(r^{-q})$, $q>n$, then Theorem \ref{scdt}
applies and we once again have $m'(t)=0$.

Combining the results in section 3, and 4, we have proved Theorem 1.

%In general, we need to assume more decay condition of the curvature
%to guarantee the vanishing of $m'(t)$.

$ \ $

Let us now make an observation. Using the divergence theorem we have
$$
m'(t)=\o{1}{4\omega_n} \int_{B_{r\rightarrow\infty}}\Delta
Rdv_{g(t)}.
$$
Comparing this with the formula
$$
\Delta R=R_t-2|Rc|^2,
$$
and
$$
\o{d}{dt} \int_{B_r} R dv_{g(t)} = \int_{B_r} R_t dv_{g(t)} -
\int_{B_r} R^2 dv_{g(t)},
$$
 we obtain that
\begin{align}
m'(t)&=\o{1}{4\omega_n} \int_{B_{r\rightarrow\infty}}(R_t-2|Rc|^2)dv_{g(t)} \nonumber \\
&=\lim_{r\to\infty} \o{1}{4\omega_n} \left[\o{d}{dt} \int_{B_r} R
dv_{g(t)} + \int_{B_r} R^2 dv_{g(t)} -2\int_{B_{r}}|Rc|^2dv_{g(t)}
\right].
\end{align}
This yields the following result.

\begin{Prop} \label{cotsc}
Under the Ricci flow for ALE metrics with (\ref{mali})
%$\int_{M}|Rc|^2dv_{g(t)}<+\infty$ for each $t\geq 0$,
we have
$$
\o{d}{dt} \int_M R dv_{g(t)}=\int_{M}\left(2|Rc|^2-R^2 \right)
dv_{g(t)},
$$
provided that $n=3$ or $\tau > n-3$.

\end{Prop}

\section{Uniform convergence}

We now turn our attention to Theorem \ref{st}. First, we introduce
the notion of {\em uniform convergence} in our context. To this end,
we now discuss the weighted Sobolev spaces $W^{k, q}_{\tau}$ and a
certain related space ${\cal M}_{\tau}$ on an ALE space $M$, see
\cite{B86,LP}.

For $q\geq 1$ and $\tau \in \mathbb R$, the weighted Lebesque space $L^q_{\tau}(M)$ consists of locally integrable functions $u$ on $M$ for which the norm
 $$
\|u\|_{q, \tau}=\left(\int_M  |r^{-\tau}u|^q r^{-n} d {\rm vol} \right)^{1/q}
 $$
is finite. For nonnegative integer $k$, the weighted Sobolev space $W^{k, q}_{\tau}(M)$ is the set of $u$ for which $|\nabla^iu| \in L^q_{\tau}(M)$ for $0\leq i \leq k$, with the norm
 $$
\|u\|_{q, k, \tau}=\sum_{i=0}^k \| \nabla^iu \|_{q, \tau}.
 $$

For $\tau > \o{n-2}{2}$, we define ${\cal M}_{\tau}$ to be the set of all $C^{\infty}$ metrics $g$ on $M$ such that, in some asymptotic coordinates,
\be \label{calm}
g_{ij} -\delta_{ij} \in W^{k, q}_{\tau}(M), \ \ \  \ \ \    R(g) \in L^1(M).
\ee
We equip ${\cal M}_{\tau}$ with the norm
\be \label{calmn}
\|g\|_{{\cal M}_{\tau}}= \|g_{ij} -\delta_{ij} \|_{W^{k, q}_{\tau}}+\| R(g) \|_{L^1}.
\ee

Now, consider the Ricci flow $(M, g(t))$ which we assume to exist for all time $0\leq t < \infty$. Furthermore, suppose that each $(M, g(t))$ is ALE of asymptotic order $\tau>\o{n-2}{2}$ and that the scalar curvature $R(t)$ of $g(t)$ is integrable on $M$ (so that the ADM mass is well defined). i.e. $g(t) \in {\cal M}_{\tau}$.

\begin{Def} \label{uc} We say that $g(t)$ {\em converges uniformly} to $g_{\infty} \in {\cal M}_{\tau}$ as $t\ra \infty$ if $g(t)$ converges to $g_{\infty}$ in ${\cal M}_{\tau}$.
\end{Def}

We now prove Theorem \ref{st}.

\begin{proof}  This follows from the argument of \cite[Lemma 9.4]{LP}.
The key here is the following identity, first observed in \cite{S,B86}.
In terms of the asymptotic coordinate,
 \be \label{domd} R(g)= \pt_j(\pt_ig_{ij}-\pt_jg_{ii}) + O(r^{-2\tau-2}),
 \ee
where the $O(r^{-2\tau-2})$ is controlled by the
$W^{1,q}_{\tau}$-norm of $g$.

Let $\eta$ be a cut off function which is identically $1$ for large
$r$ and $0$ for $r\leq 1$ and inside. Then, by Divergence Theorem
\begin{align*} m(g) & = \lim_{R \ra \infty} \frac{1}{4 \omega_n}
\int_{S_R} \eta (\partial_i g_{ij} - \pt_j g_{ii} ) * dx_j \\
 &= \int_M \left( -\eta \nabla^* \beta + \lp \beta, \nabla \eta \rp
\right) d {\rm vol},
\end{align*}
where $\beta=(\pt_ig_{ij}-\pt_jg_{ii})\pt_j$ is the mass density
vector. Theorem \ref{st} now follows from the formula above and
(\ref{domd}).
\end{proof}

%\section{Appendix}

%Finally we prove the gap result, Corollary \ref{gt}.

%Indeed, by the long time existence result of Shi \cite{Shi89}, Ricci
%flow starting at $g$ exists for all time, and converges to a flat
%metric. By \cite{BKN89} and \cite[Corollary 5.18]{Shi89}, these
%metrics $g(t)$ are all ALE. Furthermore, according to \cite[Corollary
%5.18 and Lemma 7.4]{Shi89}, $g(t)$ in fact converges uniformly.
%Thus, our Theorem \ref{st} applies and combined with Theorem
%\ref{mt} it implies that
%$$ m(g)=0. $$
%The conclusion follows immediately from the Positive Mass Theorem.

\section*{Acknowledgements} Part of the work is done while both
authors were visiting the Nankai Institute of Mathematics, Tianjin,
China. We would like to thank Nankai Institute and its Director Weiping
Zhang for the hospitality.

X.D.:\ Department of Mathematics, University of California, SANTA
BARBARA, CA 93110, USA. \ \ \  dai@math.ucsb.edu

L.M.:\ Department of Mathematical Science, Tsinghua University,
 Peking 100084, P. R. China. \ \ \ lma@math.tsinghua.edu.cn

\end{document}